\documentstyle{article}
\newcommand{\ba}{\begin{array}}
\newcommand{\ea}{\end{array}}
\newcommand{\be}{\begin{equation}}
\newcommand{\ee}{\end{equation}}
\newcommand{\ben}{\begin{enumerate}}
\newcommand{\een}{\end{enumerate}}

\newcommand{\End}{\mbox{\rm End}\, }
\newcommand{\Ldwa}{\,\stackrel{2}{\bigwedge}}
\newcommand{\Ltrzy}{\,\stackrel{3}{\bigwedge}}
\newcommand{\w}{{\!}\wedge{\!}}
\newcommand{\sca}{\eta }

\font \msb=msbm10 scaled \magstep0
\newcommand{\rtimes}{\mbox{\msb o}\,}
\newcommand{\bR}{\mbox{\msb R} }
\newcommand{\bC}{\mbox{\msb C} }
\font \msbm=msbm7 scaled \magstep0
\newcommand{\bCm}{\mbox{\msbm C} }

\font \eul=eufm10 scaled \magstep0

\newcommand{\gotG}{\mbox{\eul g}}
\newcommand{\gotH}{\mbox{\eul h}}

\newcommand{\dr}{\delta }
\newcommand{\er}{\varepsilon }
\newcommand{\lr}{\lambda }
\newcommand{\Om}{\Omega }
\def\Dr{\Delta }
\newcommand{\ov}{\overline}

\begin{document}

\begin{center}

{\bf ON BRAIDED POISSON AND QUANTUM INHOMOGENEOUS GROUPS}

\vspace{2mm}

{\sc S. Zakrzewski}

\vspace{1.3mm}

\small{\sl Department of Mathematical Methods in Physics,
University of Warsaw} \\ \small{\sl Ho\.{z}a 74, 00-682 Warsaw, Poland}
\end{center}
\vspace{3mm}

\begin{quote}
{\small
The well known incompatibility between inhomogeneous quantum
groups and the standard $q$-deformation is shown to disappear
(at least in certain cases) when admitting the quantum group to
be braided. Braided quantum $ISO(p,N\! -\! p)$ containing
$SO_q(p,N\! -\! p)$ with $|q|=1$ are constructed for
$N=2p$, $2p+1$, $2p+2$.
Their Poisson analogues (obtained first) are presented as
an introduction to the quantum case.
}
\end{quote}

\vspace{1.5mm}

\begin{center}
{\bf 1 \ \ Introduction}
\end{center}

\vspace{1mm}

It is well known \cite{PP,pspg} that the Lorentz part of any quantum
(or Poisson) Poincar\'{e} group is triangular. This is in fact a
general feature, which excludes the standard $q$-deformation
from the context of inhomogeneous quantum groups \cite{PP1}.
In order to make the standard $q$-deformation compatible
with inhomogeneous groups one has to consider some
generalization of the notion of quantum (Poisson) group,
such as, for example, a braided quantum (Poisson) group.

The notion of a braided Hopf algebra is due to S.~Majid~\cite{Ma}.
It is a natural generalization of the notion of a Hopf algebra
when we replace the usual symmetric monoidal category of vector
spaces by a braided one (the incorporation of *-structures is
more controversial --- we follow here the approach of
\cite{SLW}). A characteristic feature of this
generalization is that the comultiplication is a morphism of
algebras when the product algebra is considered with a crossed
tensor product structure rather than the ordinary one.

On the Poisson level, it means that instead of ordinary Poisson
groups $(G,\pi )$ (where $\pi $ is such a Poisson structure on
$G$ that the group multiplication is a Poisson map from the
usual product Poisson structure $\pi\oplus\pi$ on $G\times G$ to
$\pi$ on  $G$), we consider triples $(G,\pi ,\pi _{\Join})$,
where $\pi$ is a Poisson structure on $G$ and $\pi _{\Join}$ is
a bi-vector field on $G\times G$ of the cross-type (i.e. having
zero both projections on $G$) such that
\ben
\item $\pi _{12}:=\pi\oplus\pi + \pi _{\Join}$ is a Poisson structure on
$G\times G$,
\item the group multiplication is a Poisson map from $\pi _{12}$
to $\pi $.
\een
In the next section we shall construct such structures on the
inhomogeneous orthogonal groups $ISO(p,p)$, $ISO(p,p\! + \! 1)$,
$ISO(p,p\! + \! 2)$, with the homogeneous part being
non-triangular (with standard Belavin-Drinfeld $r$-matrix).

In Sect.~3, similar result is obtained for the quantum case.

\vspace{1.5mm}

\begin{center}
{\bf 2 \ \ The Poisson case}
\end{center}

\vspace{1mm}

In this section we discuss Poisson-Lie structures (possibly
braided) on inhomogeneous orthogonal groups (in particular, on
the Poincar\'{e} group). Let $V\cong \bR ^N =\bR
^{p+(N-p)}$ be equipped with the standard scalar product $\sca$
of signature $(p,N-p)$. Special linear transformations
preserving $\sca$ form the {\em homogeneous} orthogonal group
$H:=SO(p,N\!\! -\! p)\subset GL(V)$ with the Lie algebra $\gotH
:= so (p,N\!\! -\! p)\subset \End V$. The corresponding {\em
inhomogeneous} group $G=V\rtimes H$ (with Lie algebra $\gotG
=V\rtimes \gotH$) may be identified with the set of matrices
\be\label{group}
G = \left\{ \left(\ba{c|c} h & x \\ \hline 0 & 1 \ea\right)
\in \End (V\oplus \bR ) :
h\in H,\, x\in V\right\} .
\ee
For $N>2$, any multiplicative bi-vector field  $\pi$ on $G$ is
known \cite{pspg} to be of the form $\pi (g) =\pi _r
(g):=gr-rg$, where $r\in \Ldwa \gotG$. Here $r$ has three
components,
\be
r=a+b+c\;\in \; (\Ldwa V )\, \oplus (V\w \gotH )\,\; \oplus
(\Ldwa\gotH ) .
\ee
Decomposing $(V\oplus \bR )\otimes (V\oplus \bR )=
(V\otimes V)\oplus (V\otimes \bR )\oplus (\bR \otimes V)\oplus
(\bR\otimes \bR)$ (in this order), we can write tensor product
of matrices again as matrices:
\be\label{duze}
g_1g_2=
\left(\ba{cccc} h_1h_2 & h_1x_2 & x_1h_2 & x_1x_2 \\
0 & h_1 & 0 & x_1 \\
0 & 0 & h_2 & x_2 \\
0 & 0 & 0 & 1\ea\right),\qquad
r=\left(\ba{cccc} c & -b_{21} & b & a \\
0 & 0 & 0 & 0 \\
0 & 0 & 0 & 0 \\
0 & 0 & 0 & 0 \ea\right),
\ee
where the subscripts 1,2 denote the insertion place
in the tensor product. Using this, we obtain more detailed
description of the brackets defined by $\pi $,
\be
\{ g_1,g_2\} = rg_1g_2-g_1g_2r ,
\ee
as \ $\{ h_1,h_2\}  =  ch_1h_2 -h_1h_2c$, \
$\{ x_1,h_2\}  =  cx_1h_2 + bh_2 - h_1h_2b_{21}$ \ and
$\{ x_1,x_2\}  =  cx_1x_2 + bx_2-b_{21}x_1 +a -h_1h_2a$.
It follows that with any Poisson group structure on $G$ there
is associated a Poisson group structure on $H$ (with $c$ being the
$r$-matrix) and the projection from $G$ to $H$ is a Poisson map.
As shown in \cite{pspg} (see also below), $c$ must be triangular
(hence non-standard). The problem now arises if a non-triangular
$c$ can be used to construct (at least) a braided Poisson $G$.

Let us simplify the discussion to the case when $r=c$ (note that
then the inclusion $H\subset G$ is also a Poisson map). The
brackets have now the form
\be
\{ h_1,h_2\}  =  rh_1h_2 -h_1h_2r,\qquad
\{ x_1,h_2\}  =  rx_1h_2 ,\qquad
\{ x_1,x_2\}  =  rx_1x_2 .
\ee
We shall show that these brackets are not Poisson, unless $r$ is
triangular. It is convenient to check if the Jacobi identity is
satisfied in a slightly more general case:
\be\label{piw}
\{ h_1,h_2\}  =  rh_1h_2 -h_1h_2r,\qquad
\{ x_1,h_2\}  =  wx_1h_2, \qquad
\{ x_1,x_2\}  =  rx_1x_2 ,
\ee
where $w\in\gotH\otimes \gotH$. Let
$ J(f_1,f_2,f_3):=\{\{ f_1,f_2\},f_3 \}+ \{\{ f_2,f_3\},f_1
\}+\{\{ f_3,f_1\},f_2 \} $
for any functions $f_1,f_2,f_3$. It is easy to check that
\begin{eqnarray}
J(h_1,h_2,h_3) & = & [[r,r]]h_1h_2h_3-h_1h_2h_3[[r,r]] \label{hhh}\\
J(x_1,h_2,h_3) & = &
([w_{12},w_{13}]+[w_{12}+w_{13},r_{23}])x_1h_2h_3 \label{hh}\\
J(x_1,x_2,h_2) & = & ([r_{12},w_{13}+w_{23}]+
[w_{13},w_{23}])x_1x_2h_3 \label{xx}\\
J(x_1,x_2,x_3) & = & [[r,r]]x_1x_2x_3 \label{space},
\end{eqnarray}
where $[[\cdot ,\cdot ]]$ is the bracket defined by Drinfeld:
for any $\rho \in \gotH\otimes \gotH$,
$$ [[\rho ,\rho ]]:=[\rho _{12},\rho _{13}]+ [\rho _{12},\rho
_{23}] + [\rho _{13},\rho _{23}] .$$

If $w=r$, then the Jacobi identity holds provided $[[r,r]]=0$
($r$ triangular).

If $w=r+s$, where $s$ is a symmetric invariant element of
$\gotH\otimes \gotH$ and $[[w,w]]=0$ (i.e. $r$ is
real-quasitriangular),  then the Jacobi identity is satisfied,
provided (\ref{space}) is zero, i.e. the fundamental bivector
field $r_V$ on $V$ (cf.\cite{stand}) is Poisson. We shall show that
it is Poisson for almost all $N,p$, namely when
$\gotH = so(p,N\!\! -\! p)$
is absolutely simple. Indeed, in this case all invariant symmetric
2-tensors $s$ are proportional to (the Killing element)
\be\label{stil}
\widetilde{s}^{jk}_{lm} = \sca ^{jk}\sca _{lm} -\dr ^j_m\dr ^k_l ,
\ee
and all invariant elements of $\Ltrzy \gotH$ are
proportional to $\Om :=[[\widetilde{s},\widetilde{s}]]=
[\widetilde{s}_{12},\widetilde{s}_{13}]$.
{}From (\ref{stil}) we obtain
$$
\Om ^{abc}_{jkl}=
\sca ^{ab}\sca_{jl}\dr ^c_k+\sca ^{ac}\sca _{kl}\dr ^b_j +
\sca ^{bc}_{jk}\dr ^a_l -\sca^{ab}\sca_{kl}\dr^c_j
-\sca^{bc}\sca_{jl}\dr^a_k -\sca^{ac}\sca_{jk}\dr^b_l+
\dr^a_k\dr^b_l\dr^c_j - \dr ^a_l\dr ^b_j\dr ^c_k
$$
which yields $\Om ^{abc}_{jkl}x^jx^kx^l=0$. For any classical
$r$-matrix $r$ on $\gotH$, $[[r,r]]$ must be proportional to
$\Om$ and therefore (\ref{space}) is zero.

If $\gotH = so(1,3)$, all invariant symmetric 2-tensors
are complex multiples of
\be
\widetilde{s} = X_+\otimes X_-+X_-\otimes X_+ +\frac12
H\otimes H \;\;\;\;\;\;\;(\mbox{{\em complex} tensor product}).
\ee
We use here the embedding of the complex tensor product
$\gotH\otimes _{\bCm}\gotH$ into the real $\gotH\otimes
\gotH$ as described in \cite{repoi} ($X_+,X_-,H$ is the
standard complex basis of $so(1,3)\cong sl(2,\bC )$ normalized
as in \cite{repoi}; the reader should excuse the double use of
the letter $H$). One can check easily that
\be\label{ML}
\widetilde{s}=\vec{M}\cdot \vec{M}-\vec{L}\cdot\vec{L},\qquad
-i{\widetilde{s}}=\vec{M}\cdot
\vec{L}+\vec{L}\cdot\vec{M},
\ee
where $M_i:= = \er _{ijk} e_k\otimes e{^j}$, $L_i = e_0\otimes
e{^i} + e_i\otimes e{^0}$ ($i,j,k=1,2,3$) are standard
generators of $so(1,3)$ and therefore $\widetilde{s}$ coincides
with (\ref{stil}).
All invariant 3-vectors are complex multiples of
\be
\Om=[[\widetilde{s},\widetilde{s}]]= X_+\w H\w X_-
\;\;\;(\mbox{{\em complex} products; we use} \;\;
\Ltrzy _{\bCm}\gotH\;\subset \; \Ltrzy \gotH ).
\ee
Since $\Om x_1x_2x_3=0$ and $(i\Om )x_1x_2x_3\neq 0$ (Ex.~3.3 of
\cite{stand}), $r_V$ is Poisson only if $[[r,r]]$ is (real)
proportional to $\Om $. It means that if $r_-=i\lr X_+\w X_-$
(the only possibility of non-triangular $r$, up to
automorphism; the notation of \cite{repoi}), then
$[[r,r]]=[[r_-,r_-]]=\lr ^2\Om$, hence $\lr
^2$ must be real, i.e. $\lr$ real or imaginary
(cf.~\cite{stand}).

Now we turn to the question of real-quasitriangularity.
{}From Thm.~3.3 of \cite{CGR} it follows that
real-quasitriangular (not triangular) $r$-matrices exist only in
the following three cases of $so (p,N\!\! -\! p)$:

\hspace{1.3cm}  $so(p,p)$, $so(p,p+1)$ (real split cases) \ \ and
\ $so(p,p+2)$.

For $so(1,1\! +\! 2)$ in fact every $r$-matrix
is real-quasitriangular (with suitable $s$). If it is not triangular,
then, up to automorphism, $r_-=i\lr X_+\w X_-$ and $[[r,r]]=\lr
^2 \Om$, whereas
$[[s,s]]=-\lr ^2 \Om$ for $s=i\lr \widetilde{s}$, hence
$[[r+s,r+s]]=[[r,r]]+[[s,s]]=0$.

Concluding, for real-quasitriangular $r$ such that $r_V$ is
Poisson, we have a natural Poisson structure $\pi$  on $G$
defined by (\ref{piw}),
which generalizes $\pi _r$. This structure is not
multiplicative (for $s\neq 0$). It differs from the
multiplicative structure $\pi _r$ only by the following brackets:
\be
\{ h_1,h_2\}_s=0,\qquad \{ x_1,h_2\}_s :=sx_1h_2, \qquad
\{ h_1,h_2\}_s=0 .
\ee
Denoting by $\Dr $ the comultiplication: $\Dr h = hh'$, $\Dr
x=x+hx'$ (the primed functions refer to the {\em second copy} of
$G$), we obtain
$$ \{ \Dr h_1,\Dr h_2\}_s =\Dr \{ h_1,h_2\}_s,\qquad
\{ \Dr x_1,\Dr h_2\}_s =\Dr \{ x_1,h_2\}_s ,$$
but
\be
 \{ \Dr x_1,\Dr x_2\}_s - \Dr \{ x_1,x_2 \}_s = \{ \Dr x_1,\Dr
x_2\}_s = (s-Ps)x_1h_2x'_2 ,
\ee
where $P$ is the permutation in the tensor product.
It is therefore natural to look for cross-term $\{ \cdot
,\cdot \}_{\Join}$ which is nontrivial only between $x$ and $x'$.
With such an assumption, $(G,\pi ,\pi _{\Join})$ will be a
braided Poisson group if $\{\Dr x_1,\Dr x_2\}_s +\{ \Dr x_1,\Dr
x_2\}_{\Join}=0$, i.e.
\be\label{sPs}
 (s-Ps)x_1h_2x'_2 + h_2\{ x_1,x'_2\}_{\Join} +
h_1\{x'_1,x_2\}_{\Join} =0.
\ee
Consider first the generic $s$ which is proportional to
(\ref{stil}): $s=\nu\widetilde{s}$.
Since $\widetilde{s}-P\widetilde{s}=I-P$, (\ref{sPs}) is equivalent
to
\be
\nu (x_1h_2x'_2 -x_2h_1x'_1) = h_2\{ x'_2,x_1\}_{\Join} -
h_1\{x'_1,x_2\}_{\Join} ,
\ee
which is satisfied by
\be\label{Join}
\{ x'_2,x_1\}_{\Join} = \nu x_1x'_2\qquad (\mbox{more
explicitly:}\; \{ (x')^k,x^j\}_{\Join} =\nu x^j(x')^k).
\ee
One has only to check that $\pi \oplus \pi  +\pi _{\Join}$
is a Poisson bracket on $G\times G$, but this is true:
\begin{eqnarray*}
J (x_1,x_2,x'_3) & = & \{ rx_1x_2,x'_3\} + \{ x_2x'_3,x_1\}-
\{ x'_3x_1,x_2\} \\
 & = & 2r_{12}x_2x'_3+r_{21}x_2x_1x'_3
-x_2x_1x'_3+x'_3x_2x_1-r_{12}x'_3x_1x_2=0, \\
J(x_1,x'_2,h_3) & = & \{ x_1x'_2,h_3\} +\{ -w_{13}x_1h_3,x'_2\} =
 w_{13}x_1h_3x'_2 -w_{13}x_1x'_2h_3 =0
\end{eqnarray*}
(here $\{ \cdot ,\cdot \}$ denotes the full bracket on $G\times
G$ defined by $\pi \oplus \pi + \pi _{\Join}$).

In the Lorentz case $\gotH =so(1,3)$, apart from the generic
case $s=\nu \widetilde{s}$, one has to consider also the case when
$s=\nu i\widetilde{s}$. Using formula (\ref{ML})
for $i\widetilde{s}$, it is easy to see that $i\widetilde{s}
-Pi\widetilde{s} = 2i\widetilde{s}$ and (\ref{sPs}) has no
solutions. Thus the case of real $\lr $ in $r_-=i\lr X_+\w X_-$,
which corresponds to real $q$ in the quantum case (in
particular, quantum double of $SU_q(2)$), is excluded.
It means that from the list of $r$-matrices on $so(1,3)$ in
\cite{repoi}, only combinations of $(X_+\w X_- -JX_+\w JX_-)$
and $JH\w H$ fall in our scheme.

Finally, it is interesting to note that
\ben
\item
the one-parameter group of automorphisms of $G$ (dilations),
$$ t(h,x):= (h,e^tx)\qquad\mbox{for}\;\; t\in \bR,$$
preserves $\pi$ (because (\ref{piw}) is homogeneous in $x$),
\item the braiding bivector field $\pi _{\Join}$ described by
(\ref{Join}) is nothing else but the antisymmetrization of the
fundamental tensor field on $G\times G$ obtained by the action
of the real-quasitriangular element
$$ \nu e_1\otimes e_1\in \bR\otimes \bR\qquad (e_1\;\mbox{is the
basic vector of}\;\;\bR ).$$
\een
Similar property is satisfied by the {\em cobracket} $\dr$ on $\gotG$,
obtained by linearization of $\pi $ at the group unit. It follows
that $(\gotG , \dr )$ is an example of a {\em braided-Lie
bialgebra} \cite{M2} (in the category of modules over
quasitriangular $\bR $). $(G,\pi )$ will certainly be an example
of a braided Poisson-Lie group, when  the theory presented in
\cite{M2} will be extended from Lie algebras to Lie groups.

\vspace{1.5mm}

\begin{center}
{\bf 3 \ \ The quantum case}
\end{center}

\vspace{1mm}

Real-(co)quasitriangular quantum $SO(p,p)$ and $SO(p,p+1)$ are
introduced in \cite{FRT} and $SO(p,p+2)$ in \cite{it}. They all
can be described by relations of the form
\be\label{homor}
Wh_1h_2=h_2h_1W, \qquad h_1h_2 \eta  =
\eta,\;\;  \eta ' h_1h_2=\eta ',\;\;\; h=h^*,
\ee
where
\be\label{R-matrix}
\hat{W}=PW = qP^{(+)}-q^{-1}P^{(-)}+q^{1-N}P^{(0)}
\ee
is the standard $R$-matrix for the orthogonal series (here
$P^{(+)}$, $P^{(-)}$ and $P^{(0)}$ are the spectral projections
corresponding to symmetric (traceless), antisymmetric and
proportional to the metric elements of $V\otimes V$)
with $|q|=1$ and $\eta '$ ($\eta $) is a deformed covariant
(contravariant) metric. For
$q=1+i\er +\ldots$ we have $W=I+i\er w +\ldots $, where $w$
satisfies the classical Yang Baxter equation. To the
skew-symmetric classical $r$-matrix $r=(w-w_{21})/2$ there
corresponds the involutive intertwiner
$$ \hat{R}:= I-2P^{(-)},\qquad R=P\hat{R}=I+i\er r+\ldots $$
(note that $R$ can be used instead of $W$ in (\ref{homor})).

Passing to the inhomogeneous group (\ref{group}), we expect that the
commutation relations for $g$ should be
\be
{\cal R}g_1g_2=g_2g_1{\cal R},\qquad \mbox{where}\;\;\;
{\cal R}=\left(\ba{cccc} R & 0 & 0 & 0 \\
0 & I & 0 & 0 \\
0 & 0 & I & 0 \\
0 & 0 & 0 & 1 \ea\right)
\ee
(this corresponds to $r$ given in (\ref{duze}) when $a=0$, $b=0$).
Using the form of $g_1g_2$ as in (\ref{duze}), we obtain
\be\label{triang}
Rh_1h_2=h_2h_1R,\qquad x_2h_1=Rh_1x_2,\qquad
h_2x_1=Rx_1h_2,\qquad x_2x_1=Rx_1x_2.
\ee
The two equalities in the middle are equivalent, due to the
involutivity of $\hat{R}$. The last equality provides defining
relations for the {\em quantum orthogonal vector space}
\cite{FRT,it}. These relations are consistent: the corresponding
algebra of polynomials has the classical size. Also the first
equality gives consistent relations in this sense. It remains to
check the consistency of the `cross-relations' with other ones. From
\begin{eqnarray}
R_{12}R_{13}R_{23}h_1h_2x_3 & = & x_3h_2h_1R_{12} =
R_{23}R_{13}R_{12}h_1h_2x_3, \label{R1}\\
R_{12}R_{13}R_{23}h_1x_2x_3 & = & x_3x_2h_1 =
R_{23}R_{13}R_{12}h_1x_2x_3, \label{R2}
\end{eqnarray}
it follows that $R$ should satisfy the Yang Baxter equation,
hence $q=1$ (the triangular case).
As in the Poisson case, we postulate then a modification
of (\ref{triang}) as follows:
\be\label{W}
Rh_1h_2=h_2h_1R,\qquad x_2h_1=W'h_1x_2,\qquad
 x_2x_1=Rx_1x_2, \;\;\;\;
\ee
with some matrix $W'$. Instead of (\ref{R1})--(\ref{R2}), we
have now
\begin{eqnarray*}
R_{12}W'_{13}W'_{23}h_1h_2x_3 & = & x_3h_2h_1R_{12} =
W'_{23}W'_{13}R_{12}h_1h_2x_3, \\
W'_{12}W'_{13}R_{23}h_1x_2x_3 & = & x_3x_2h_1 =
R_{23}W'_{13}W'_{12}h_1x_2x_3.
\end{eqnarray*}
For the consistency of different ways of ordering, we postulate that
\be\label{post}
W'_{12}W'_{13}W'_{23} = W'_{23}W'_{13}W'_{12}\qquad\mbox{and}\;
\hat{R}\;(\mbox{or}\;P^{(-)})\;\mbox{is a function of}\;\; \hat{W'}=PW'.
\ee
This is fulfilled if $\hat{W}'$ a scalar multiple of $\hat{W}$
(it is also possible that $\hat{W}'$ is a scalar multiple of
$\hat{W}^{-1}$; this corresponds to the change $s\mapsto -s$ in
the Poisson case). The scalar coefficient is not arbitrary, due
to the following two conditions:
\ben
\item From the reality requirement ($h^*=h$, $x^*=x$) it follows that
$x_2h_1=W'h_1x_2 $ implies $h_1x_2=\ov{W'}x_2h_1$, hence
$x_2h_1=W'\ov{W'}x_2h_1$ and we have to assume that
\be\label{reality}
W'\ov{W'}=I.
\ee
\item
Since $x_3\eta  _{12}=x_3h_1h_2\eta _{12}=W'_{13}W'_{23}h_1h_2x_3\eta
_{12}= W'_{13}W'_{23}\eta _{12}x_3$, we have also the following
condition of compatibility of $W'$ with the metric:
\be\label{compa}
W'_{13}W'_{23}\eta _{12}=\eta _{12}.
\ee
\een
Both conditions are satisfied by $W'=W$ (another solution,
$W'=-W$, has no proper classical limit). The first condition
follows from
$$ \ov{W(q)}=W(\ov{q})=W(q^{-1})=W(q)^{-1}$$
(cf. \cite{FRT}; recall that $|q|=1$). The second coincides with
formula (2.21) in \cite{Po}. Thus, in the sequel we set $W'=W$.

It is easy to see that the comultiplication preserves first two
relations in (\ref{W}), for instance $\Dr x_2\Dr h_1$ equals
$$(x_2+h_2x'_2)h_1h_1'=Wh_1x_2h'_1+h_2h_1Wh'_1x'_2=
Wh_1h'_1x_2+Wh_1h_2h'_1x'_2=W\Dr h_1\Dr x_2.$$
This will be true also for a nontrivial braiding of the type
\be
x'_2x_1=Bx_1x'_2,
\ee
which on the other hand may be used to remove the inconsistency
related to the preservation of the third relation:
 $P^{(-)}x_1x_2=0$. We shall find now the
condition under which $P^{(-)}\Dr x_1\Dr x_2=0$. The first
two terms in
$$
\Dr x_1\Dr x_2 = (x_1+h_1x'_1)(x_2+h_2x'_2)=x_1x_2
+h_1x'_1h_2x'_2 + x_1h_2x'_2 + h_1x'_1x_2
$$
are annihilated by $P^{(-)}$ (second, because
$P^{(-)}h_1h_2x'_1x'_2= h_1h_2P^{(-)}x'_1x'_2=0$).
The sum of the last two terms is equal
$$
(\hat{W}h_1x_2x'_1 + h_1x'_1x_2)^{jk}=\hat{W}^{jk}_{ab}h^a_cx^bx'^c+
h^j_lB^{kl}_{bc}x^bx'^c=
(\hat{W}^{jk}_{ab}\dr ^l_c+\dr ^j_aB^{kl}_{bc})h^a_lx^bx'^c,
$$
hence our condition is
\be
P^{(-)}_{12}(\hat{W}_{12}+B_{23})=0.
\ee
If
\be\label{main}
P^{(-)}(\hat{W}+\sigma I)=0\qquad\mbox{for some}\;\;\sigma ,
\ee
then $B=\sigma I$ is a solution of our problem and the non-trivial
cross-relations are the following: $x'^jx^k=\sigma
x^kx'^j$. We call (\ref{main}) the {\em spectral condition}.
Taking into account that $P^{(-)}$ is a projection and a
function of $\hat{W}$, it means that $P^{(-)}$ is a spectral
projection of $\hat{W}$ corresponding to a single eigenvalue.
This is of course satisfied for (\ref{R-matrix}), with $\sigma
=q^{-1}$. 

We conclude that relations (\ref{W}) with $W'=W$ and braiding
\be
x'^jx^k=q^{-1}x^kx'^j
\ee
define a braided quantum $ISO(p,N\! -\! p)$, which contains
$SO_q(p,N\! -\! p)$.

\end{document}